\documentstyle[11pt]{article}
\begin{document}
\author{{Shiqiu Zheng$^{1}$\thanks{E-mail: shiqiumath@163.com (S. Zheng).}\ , \ \ Gaofeng Zong$^{2}$\thanks{E-mail: gf\b{ }zong@126.com (G. Zong).}}\\
  \\
\small(1, School of Science, Tianjin University of Science and Technology, Tianjin 300457, China)\\
\small(2, School of Mathematics and Quantitative Economics, \\ \small Shandong University of Finance and Economics, Jinan 250014, China)
}
\date{}
\title{\textbf{BSDEs and SDEs with time-advanced and -delayed coefficients}\thanks{The first author is supported by National Natural Science Foundation of China (No. 11571024) and a program of Hebei province (No. QN2017116). The second author is supported by the National Natural Science Foundation of China (No. 11501325).}}\maketitle

\textbf{Abstract:}\quad This paper introduces a class of backward stochastic differential equations (BSDEs), whose coefficients not only depend on the value of its solutions of the present but also the past and the future. For a sufficiently small time delay or a sufficiently small Lipschitz constant, the existence and uniqueness of such BSDEs is obtained. As an adjoint process, a class of stochastic differential equations (SDEs) is introduced, whose coefficients also depend on the present, the past and the future of its solutions. The existence and uniqueness of such SDEs is proved for a sufficiently small time advance or a sufficiently small Lipschitz constant. A duality between such BSDEs and SDEs is established.\\



\section{Introduction}
Peng and Yang [9] introduced a type of backward stochastic differential equations (BSDEs), called anticipated BSDEs, whose coefficients depend not only on the values of solutions of the present but also the future. Delong and Imkeller [6] introduced a type of BSDEs whose coefficients depend not only on the values of solutions of the present but also the past. The two types of BSDEs have been applied in many problems arising from finance and stochastic control. It is natural to consider the BSDEs with time-advanced and -delayed coefficients, i.e. their coefficients not only depend on the value of its solutions of the present but also the past and the future. In fact, Cheridito and Nam [2] had studied BSDEs whose coefficients can depend on the whole path of the solutions in a very general case. In this paper, we study BSDEs with time-advanced and -delayed coefficients in a more specific case. We consider the case that the coefficient $g(t,\cdot,\cdot)$ depends on the solution on the interval $[t-l,t+u]$, where $l\geq0$ and $u\geq0$ are the time-delayed parameter and the time-advanced parameter, respectively. For a sufficiently small time delay $l$ or a sufficiently small Lipschitz constant, the existence and uniqueness of solution of such BSDEs is obtained. Since we use an estimate different from [6], our result is independent of the terminal time $T,$ which is different from the corresponding results in [6]. More generally, we consider the BSDE whose coefficient $g(t,\cdot,\cdot)$ depends on the solution on the interval $[-l,T+u]$, and prove such BSDEs have a unique solution for a sufficiently small $T+l$ or a sufficiently small Lipschitz constant.

Recently, to study stochastic maximum principle, Chen and Huang [1] introduces a type of SDEs with time-advanced coefficients, whose coefficients depend on the future of its solutions. As an adjoint process of BSDEs with time-advanced and -delayed coefficients, in this paper, we introduce a type of SDEs, whose coefficients not only depend on the value of its solutions of the present but also the past and the future. We consider the case that the coefficients ${b}(t,\cdot)$ and ${\sigma}(t,\cdot)$ depend on the solution on the interval $[t-l,t+u]$, where $l\geq0$ and $u\geq0$ are the time-delayed parameter and the time-advanced parameter, respectively. For a sufficiently small time advance $u$ or a sufficiently small Lipschitz constant, the existence and uniqueness of solution of such SDEs is proved. More generally, we consider the SDE whose the coefficients ${b}(t,\cdot)$ and ${\sigma}(t,\cdot)$ depend on the solution on the interval $[t_0-l,T+u]$, and prove such SDE has a unique solution for a sufficiently small $T+u-t_0$ or a sufficiently small Lipschitz constant.

In our opinion, the BSDEs and SDEs with time-advanced and -delayed coefficients can be well applied in the dynamical problems depending on the past, the present and the future. Exploring their applications will be our future work. In fact, in the determinate case, d'Albis et al. [4][5] had introduced the mixed-type functional differential equations (coefficients depend on the past, the present and the future), and given their applications in economic problems. Our proofs of the existence and uniqueness of such BSDEs and SDEs both use the contractive mapping method based on the two estimates for It\^{o}'s type processes, respectively, which are mainly inspired by Peng [8]. We also obtain the continuous dependence properties and comparison theorems for such BSDEs and SDEs, and a duality between them. From this paper, one can see an interesting similarity between the results and arguments of such BSDEs and SDEs.
\section{BSDEs with time-advanced and -delayed coefficients}
In this paper, we consider a complete probability space  $(\Omega,\cal{F},\mathit{P})$ on which a $d$-dimensional standard Brownian motion ${{(B_t)}_{t\geq
0}}$ is defined. Let $({\cal{F}}_t)_{t\geq 0}$ denote
the natural filtration generated by ${{(B_t)}_{t\geq 0}}$, augmented
by the $\mathit{P}$-null sets of ${\cal{F}}$. Let $|z|$ denote its
Euclidean norm, for $\mathit{z}\in {\mathbf{R}}^d.$ Let $T>0, l\geq0, u\geq0$ be given constants. For $s\in[0,l]$ and $t\in[0,T+u],$ we define the following usual spaces:

${\cal{L}}^2([-s,t];{\mathbf{R}}^d)=\{f(t):$ Lebesgue measurable ${\mathbf{R}}^d$-valued function defined on $[-s,t]$ and $\int_{-s}^t|f(r)|^2dr<\infty$$\};$

$L^2({\mathcal {F}}_t;{\mathbf{R}}^d)=\{\xi: {\cal{F}}_t$-measurable
${\mathbf{R}}^d$-valued random variable such that $E\left[|\xi|^2\right]<\infty\};$

$L^2_{\cal{F}}(0,t;{\mathbf{R}}^d)=\{\psi: {\mathbf{R}}^d$-valued progressively measurable
process such that $E\int_0^t|\psi_r|^2dr
<\infty \};$

${\cal{S}}^2_{\cal{F}}(0,t;{\mathbf{R}}^d)=\{\psi:$ continuous
process in $L^2_{\cal{F}}(0,t;{\mathbf{R}}^d)$ such that $E[\sup_{t\in[0,T]}|\psi_t|^2]<\infty\};$

$\overline{L}^2_{\cal{F}}(-s,t;{\mathbf{R}}^d)=\{\psi:\ (\psi_r)_{r\in[0,t]}\in L^2_{\cal{F}}(0,t;{\mathbf{R}}^d) \ \textmd{and}\ (\psi_r)_{r\in[-s,0]}\in{\cal{L}}^2([-s,0];{\mathbf{R}}^d)\},\ s\geq0.$

Note that for convenience, we denote $\{X_{r}\}_{-l\leq r\leq T+u}$ by $\textbf{X}$ and when $d=1,$ denote $L^2({\mathcal {F}}_t;{\mathbf{R}}^d)$ by $L^2({\mathcal {F}}_t)$ and make the same treatment for above notations of other spaces. Clearly, $\overline{L}^2_{\cal{F}}(-s,t;{\mathbf{R}}^d)$ is a Banach spaces, and will be considered as the spaces of solutions of BSDEs and SDEs considered in this paper.

Now, we consider a function $g$
$${g}(\omega,t,y,z): \Omega\times[0,T]\times\overline{L}^2_{\cal{F}}(-l,T+u)\times \overline{L}^2_{\cal{F}}(-l,T+u;{\mathbf{R}}^d)\longmapsto {\mathbf{R}},$$  such that for each $(\textbf{Y},\textbf{Z})\in \overline{L}^2_{\cal{F}}(-l,T+u)\times \overline{L}^2_{\cal{F}}(-l,T+u;{\mathbf{R}}^d),$
$({g}(t,\textbf{Y},\textbf{Z}))_{t\in [0,T]}$ is a progressively measurable process. For the function $g$, we make the following assumptions:
\begin{itemize}
  \item (A1). There exist a constant $K>0$ and a probability measure $\lambda$ defined on $[-l,u]$, such that for each $t\in[0,T]$ and each $(\textbf{y},\textbf{z}), (\tilde{\textbf{y}},\tilde{\textbf{z}})\in \overline{L}^2_{\cal{F}}(-l,T+u)\times \overline{L}^2_{\cal{F}}(-l,T+u;{\mathbf{R}}^d),$
  \begin{eqnarray*}
   |{g}(t,\textbf{y},\textbf{z})-{g}(t,\tilde{\textbf{y}},\tilde{\textbf{z}})|
   \leq KE\left[\int_{-l}^{u}(|y_{t+r}
    -\tilde{y}_{t+r}|+|z_{t+r}
    -\tilde{z}_{t+r}|)\lambda(dr)|{\cal{F}}_t\right];
\end{eqnarray*}
    \item (A2). There exists a constant $K_1>0$, such that for each $(\textbf{y},\textbf{z}), (\tilde{\textbf{y}},\tilde{\textbf{z}})\in \overline{L}^2_{\cal{F}}(-l,T+u)\times \overline{L}^2_{\cal{F}}(-l,T+u;{\mathbf{R}}^d),$
  \begin{eqnarray*}
   E\left[\int_0^T|{g}(t,\textbf{y},\textbf{z})-{g}(t,\tilde{\textbf{y}},\tilde{\textbf{z}}) |^2dt\right]
   \leq K_1E\left[\int_{-l}^
   {T+u}|y_{t}
    -\tilde{y}_{t}|^2+|z_{t}
    -\tilde{z}_{t}|^2dt\right];
\end{eqnarray*}
  \item (A3). $g(t,\textbf{0},\textbf{0})\in L^2_{\cal{F}}(0,T).$
\end{itemize}
\

Now, we give a remark on the above assumptions.\\\\
\textbf{Remark 1}  $(i)$\ \ We consider two special cases of (A1). If for each $t\in[0,T]$ and  each $(\textbf{y},\textbf{z})\in \overline{L}^2_{\cal{F}}(-l,T+u)\times \overline{L}^2_{\cal{F}}(-l,T+u;{\mathbf{R}}^d),$ $g(t,\textbf{y},\textbf{z})$ is dependent only on $(t,\textbf{y},z_t)$, then we denote $g(t,\textbf{y},\textbf{z})$ by ${g}(t,\textbf{y},z_t)$. We give this Lipschitz condition (A1)' for it.
\begin{itemize}
 \item (A1)'. There exist a constant $K'>0$ and a probability measure $\lambda'$ defined on $[-l,u]$, such that for each $t\in[0,T]$ and each $(\textbf{y},z_t), (\tilde{{\textbf{y}}}_t,\tilde{{z}}_t)\in \overline{L}^2_{\cal{F}}(-l,T+u)\times L^2({\mathcal {F}}_t;{\mathbf{R}}^d),$
  \begin{eqnarray*}
   |{g}(t,\textbf{y},z_t)-{g}(t,\tilde{{\textbf{y}}},\tilde{{z}}_t)|
   \leq K'\left(E\left[\int_{-l}^{u}|y_{t+r}
    -\tilde{y}_{t+r}|\lambda'(dr)|{\cal{F}}_t\right]+|z_{t}
    -\tilde{z}_{t}|\right).
  \end{eqnarray*}
\end{itemize}
If for each $t\in[0,T]$ and  each $(\textbf{y},\textbf{z})\in \overline{L}^2_{\cal{F}}(-l,T+u)\times \overline{L}^2_{\cal{F}}(-l,T+u;{\mathbf{R}}^d),$ $g(t,\textbf{y},\textbf{z})$ is dependent only on $(t,y_t,z_t)$, then we denote $g(t,\textbf{y},\textbf{z})$ by ${g}(t,y_t,z_t)$. We give this standard Lipschitz condition (A1)'' (see Peng [8]) for it.
\begin{itemize}
 \item (A1)''. There exist a constant $K''>0$, such that for each $t\in[0,T]$ and each $(y_t,z_t), (\tilde{{y}}_t,\tilde{{z}}_t)\in {L}^2({\cal{F}}_t)\times L^2({\mathcal {F}}_t;{\mathbf{R}}^d),$
  \begin{eqnarray*}
   |{g}(t,y_t,z_t)-{g}(t,\tilde{{y}}_t,\tilde{{z}}_t)|
   \leq K''(|y_{t}
    -\tilde{y}_{t}|+|z_{t}
    -\tilde{z}_{t}|).
  \end{eqnarray*}
\end{itemize}
We can check that (A1)''implies (A1)', and (A1)' implies (A1).

 $(ii)$\ \ (A1) implies (A2). In fact, if $g$ satisfies (A1), then by (A1) and Fubini's theorem, for each $(\textbf{y},\textbf{z}), (\tilde{\textbf{y}},\tilde{\textbf{z}})\in \overline{L}^2_{\cal{F}}(-l,T+u)\times \overline{L}^2_{\cal{F}}(-l,T+u;{\mathbf{R}}^d),$ we have
\begin{eqnarray*}
 E\left[\int_0^T|{g}(t,\textbf{y},\textbf{z})-{g}(t,\tilde{\textbf{y}},\tilde{\textbf{z}}) |^2dt\right]&\leq&2K^2E\left[\int_0^T \int_{-l}^{u}(|y_{t+r}-\tilde{y}_{t+r}|^2+|z_{t+r}-\tilde{z}_{t+r}|^2)\lambda(dr) dt\right]\\&=&2K^2 E\left[\int_{-l}^{u}\int_0^T(|y_{t+r}-\tilde{y}_{t+r}|^2+|z_{t+r}-\tilde{z}_{t+r}|^2)ds\lambda(dr)\right]\\&=&2K^2 E\left[\int_{-l}^{u}\int_r^{T+r}(|y_{v}-\tilde{y}_{v}|^2+|z_{v}-\tilde{z}_{v}|^2
)dv\lambda(dr)\right]\\&\leq& 2K^2E\left[\int_{-l}^{T+u}(|y_{v}-\tilde{y}_{v}|^2+|z_{v}-\tilde{z}_{v}|^2)dv\right].
\end{eqnarray*}
Thus, $g$ satisfies (A2).

$(iii)$ For $t\in[0,T],\ (\textbf{y},\textbf{z})\in \overline{L}^2_{\cal{F}}(-l,T+u)\times \overline{L}^2_{\cal{F}}(-l,T+u;{\mathbf{R}}^d),$ let
\begin{eqnarray*}
&&g_1(t,\textbf{y},\textbf{z})=E\left[|y_{t+\delta_1}|+z_{t+\delta_2}+\int_{\delta_3}^{\delta_4}(y_{t+r}+{z}_{t+r})dr|{\cal{F}}_t\right],\ \ \delta_1,\delta_2,\delta_3,\delta_4\in[-l,u];\\
&&g_2(t,\textbf{y},\textbf{z})=E\left[y_{t+\delta_1}+|z_{t+\delta_2}|+\int_{\delta_3}^{\delta_4}(y_{r}+{z}_{r})dr|{\cal{F}}_t\right],\ \ \delta_1,\delta_2\in[-l,u],\delta_3,\delta_4\in[-l,T+u].
\end{eqnarray*}
 By a predictable projection theorem (see Jacod and Shiryaev [7, page 23]), one can see $g_1$ and $g_2$ are progressively measurable. Moreover, we can check that $g_1$ satisfy (A1) and (A3), and $g_2$ satisfies (A2) and (A3).\\

For convenience, if $(\xi_t,\eta_t)\in L^2_{\cal{F}}(T,T+u)\times L^2_{\cal{F}}(T,T+u;{\mathbf{R}}^d)$ with $\xi_T\in L^2({\mathcal {F}}_T)$, then we call $(\xi_t,\eta_t)$ satisfies terminal condition. Let $(\xi_t,\eta_t)$ satisfy terminal condition, we consider the following BSDE:
$$\left\{
    \begin{array}{ll}
      Y_t=\xi_T +\int_t^Tg(s,\textbf{Y},\textbf{Z})
             ds-\int_t^TZ_sdB_s,\ \ \ t\in[0,T);\\
      (Y_t, Z_t)=(\xi_t, \eta_t),\ t\in[T,T+u]\ \ \textmd{and} \ \
      (Y_t, Z_t)=(Y_0, 0), \ t\in[-l,0).
    \end{array}
  \right.\eqno(1)$$
The solution of BSDE (1) is a pair $(\textbf{Y},\textbf{Z})$ satisfying (1) and $(Y_t, Z_t)_{t\in[0,T]}\in {\cal{S}}^2_{\cal{F}}(0,T)\times L^2_{\cal{F}}(0,T;{\mathbf{R}}^d)$, which depends only on the parameter $(g,T, \xi, \eta, l,u).$ Clearly, the coefficient $g$ not only depends on the value of its solutions of the present but also the past and the future.\\

From Peng [8, Lemma 3.1], we can get the following Lemma 2.1.\\\\
\textbf{Lemma 2.1}\textit{\ Let $(\xi_t,\eta_t)$ satisfy terminal condition and $g_0(s)\in {L}^{2}_{\cal{F}}(0,T).$ Then the BSDE (1) with coefficient $g_0(s)$
has a unique solution $(\textbf{Y},\textbf{Z}),$ and the following estimate}
$$|Y_0|^2e^{-\beta l}+E\left[\int_{-l}^T(\frac{\beta}{2}|Y_s|^2+|Z_s|^2)e^{\beta s} ds\right]\leq E\left[|\xi|^2e^{\beta T}\right]+\frac{2}{\beta}E\left[\int_{0}^T|g_0(s)|^2e^{\beta s}ds\right],\eqno(2)$$
\emph{holds true for an arbitrary constant $\beta>0$. We also have}
$$E\left[\sup_{0\leq t \leq T}|Y_t|^2\right]\leq CE\left[|\xi_T|^2+\int_0^T|g_0(s)|^2ds\right],\eqno(3)$$
\textit{where $C>0$ is a constant depending only on $T$.}\\\\
\textit{Proof.} By Peng [8, Lemma 3.1],  such BSDE has a unique solution, and for an arbitrary constant $\beta>0$ and $t\in[0,T],$ we have
$$|Y_t|^2+E\left[\int_{t}^T(\frac{\beta}{2}|Y_s|^2+|Z_s|^2)e^{\beta (s-t)} ds|{\cal{F}}_t\right]\leq E\left[|\xi|^2e^{\beta T}|{\cal{F}}_t\right]+\frac{2}{\beta}E\left[\int_{t}^T|g_0(s)|^2e^{\beta (s-t)}ds|{\cal{F}}_t\right].\eqno(4)$$
Since $Y_s=Y_0,\ s\in[-l,0),$ we have
$$|Y_0|^2=|Y_0|^2e^{-\beta l}+\beta\int_{-l}^0|Y_s|^2e^{\beta s} ds.\eqno(5)$$
By the fact $Z_s=0,\ s\in[-l,0)$, (4) and (5), we can get (2).  By BDG inequality and (2), we can get (3). This proof is complete. $\Box$\\

The following is the main result of this section. It gives two sufficient conditions, under which BSDE (1) has a unique solution. \\\\
\textbf{Theorem 2.2}\textit{\ Let $(\xi_t,\eta_t)$ satisfy terminal condition. Then we have}

\textit{(i) If $g$ satisfies (A1) and (A3), and there exists a constant $\beta\geq2$ such that $\frac{4 K^2e^{\beta l}}{\beta}<1,$ then BSDE (1)
has a unique solution, where $K$ and $l$ are the constants in (A1).}

\textit{(ii) If $g$ satisfies (A2) and (A3), and there exists a constant $\beta\geq2$ such that $\frac{2K_1e^{\beta (T+l)}}{\beta}<1,$  then BSDE (1)
has a unique solution, where $K_1$ and $l$ are the constants in (A2).}\\\\
\emph{Proof.} We will prove this theorem using a contractive method. We can check that $\overline{L}^2_{\cal{F}}(-l,T+u)\times \overline{L}^{2}_{\cal{F}}(-l,T+u;{\mathbf{R}}^d)$ is a Banach space with the norm:
$$\|(\cdot,\cdot)\|_\beta=E\left[\int_{-l}^{T+u}(|\cdot|^2+|\cdot|^2)e^{\beta s} ds\right]^{\frac{1}{2}},$$
where $\beta>0$ is a constant. If $g$ satisfies (A2) and (A3), then for each $(\textbf{y},\textbf{z})\in \overline{L}^2_{\cal{F}}(-l,T+u)\times \overline{L}^{2}_{\cal{F}}(-l,T+u;{\mathbf{R}}^d),$ we have
\begin{eqnarray*}
E\left[\int_{0}^T|g(s,\textbf{y},\textbf{z})|^2ds\right]&\leq& 2E\left[\int_{0}^T|g(s,\textbf{0},\textbf{0})|^2ds\right]+2K_1E\left[\int_{-l}^{T+u}(|{y}_{s}|^2+|{z}_{s}|^2)ds\right]\\
&<&\infty.
\end{eqnarray*}
If $g$ satisfies (A1) and (A3), then by (ii) in Remark 1, we also have $$E\left[\int_{0}^T|g(s,\textbf{y},\textbf{z})|^2ds\right]<\infty.$$
Thus if $g$ satisfies (A1) (or (A2)) and (A3), by Lemma 2.1, we can define a mapping $\phi$ from
$\overline{L}^2_{\cal{F}}(-l,T+u)\times \overline{L}^{2}_{\cal{F}}(-l,T+u;{\mathbf{R}}^d)$
into itself by setting $(\textbf{Y},\textbf{Z}):=\phi(\textbf{y},\textbf{z}),$
where $(\textbf{Y},\textbf{Z})$ is the solution of the BSDE
$$\left\{
    \begin{array}{ll}
      Y_t=\xi_T +\int_t^Tg(s,\textbf{y},\textbf{z})
             ds-\int_t^TZ_sdB_s,\ \ \ t\in[0,T);\\
      (Y_t, Z_t)=(\xi_t, \eta_t),\ t\in[T,T+u]\ \ \textmd{and} \ \
      (Y_t, Z_t)=(Y_0, 0), \ t\in[-l,0).
    \end{array}
  \right.$$
For any $(\textbf{y}^1,\textbf{z}^1), (\textbf{y}^2,\textbf{z}^2)\in\overline{L}^2_{\cal{F}}(-l,T+u)\times \overline{L}^{2}_{\cal{F}}(-l,T+u;{\mathbf{R}}^d),$ let $(\textbf{Y}^i,\textbf{Z}^i):=\phi(\textbf{y}^i,\textbf{z}^i), i=1,2.$ We set
$(\hat{Y}_t,\hat{Z}_t):=({Y}^1_t-{Y}^2_t,{Z}^1_t-{Z}^2_t)$ and $(\hat{y}_t,\hat{z}_t):=({y}^1_t-{y}^2_t,{z}^1_t-{z}^2_t).$

\textbf{Proof of (i)}:

By Lemma 2.1, (A1) and Fubini's theorem, we can deduce, for an arbitrary constant $\beta>0,$
\begin{eqnarray*}
\ \ \ \ \  E\left[\int_{-l}^{T+u}(\frac{\beta}{2}|\hat{Y}_s|^2+|\hat{Z}_s|^2)e^{\beta s} ds\right]&\leq& \frac{2}{\beta}E\left[\int_0^T|g(s,{\textbf{y}}^1,{\textbf{z}}^1)-g(s,{\textbf{y}}^2,{\textbf{z}}^2\})|^2e^{\beta s}ds\right]\\&\leq& \frac{2}{\beta}E\left[2K^2\int_0^T\int_{-l}^{u}(|\hat{y}_{s+r}|^2+|\hat{z}_{s+r}|^2)\lambda(dr) e^{\beta s}ds\right]\\&=&\frac{2}{\beta}E\left[2K^2\int_{-l}^{u}\int_0^T(|\hat{y}_{s+r}|^2+|\hat{z}_{s+r}|^2) e^{\beta s}ds\lambda(dr)\right]\\&=&\frac{2}{\beta}E\left[2K^2\int_{-l}^{u}e^{-\beta r}\int_r^{T+r}(|\hat{y}_v|^2
+|\hat{z}_v|^2)e^{\beta v}dv\lambda(dr)\right]\\&\leq& \frac{4K^2}{\beta}e^{\beta l}E\left[\int_{-l}^{T+u}(|\hat{y}_v|^2
+|\hat{z}_v|^2)e^{\beta v}dv\right].\ \ \ \ \ \ \ \ \ \ \ \ \ \ \ \ \ \ (6)
\end{eqnarray*}
Thus, if there exists a constant $\beta\geq2$ such that $\frac{4 K^2e^{\beta l}}{\beta}<1,$ then there exists a constant $0<\gamma<1,$ such that
$$E\left[\int_{-l}^{T+u}(|\hat{Y}_s|^2+|\hat{Z}_s|^2)e^{\beta s} ds\right]\leq \gamma E\left[\int_{-l}^{T+u}(|\hat{y}_s|^2
+|\hat{z}_s|^2)e^{\beta s}ds\right].$$ Thus, $\|(\hat{Y}_s,\hat{Z}_s)\|_\beta<\sqrt{\gamma}\|(\hat{y}_s,\hat{z}_s)\|_\beta.$ Then by contraction mapping principle, we can obtain (i). Moreover, by Lemma 2.1, we have $(Y_t)_{t\in[0,T]}\in{\mathcal{S}}^2_{\cal{F}}(0,T)$.

\textbf{Proof of (ii)}:

By Lemma 2.1 and (A2), we can deduce, for an arbitrary constant $\beta>0,$
\begin{eqnarray*}
\ \ \ \ \  E\left[\int_{-l}^{T+u}(\frac{\beta}{2}|\hat{Y}_s|^2+|\hat{Z}_s|^2)e^{\beta s} ds\right]&\leq& \frac{2}{\beta}E\left[\int_0^T|g(s,{\textbf{y}}^1,{\textbf{z}}^1)-g(s,{\textbf{y}}^2,{\textbf{z}}^2\})|^2e^{\beta s}ds\right]\\
&\leq& \frac{2}{\beta}e^{\beta T}E\left[\int_0^T|g(s,{\textbf{y}}^1,{\textbf{z}}^1)-g(s,{\textbf{y}}^2,{\textbf{z}}^2\})|^2ds\right]
\\&\leq& \frac{2 K_1}{\beta}e^{\beta T}E\left[\int_{-l}^{T+u}(|\hat{y}_{s}|^2+|\hat{z}_{s}|^2)ds\right]
\\&\leq& \frac{2 K_1}{\beta}e^{\beta (T+l)}E\left[\int_{-l}^{T+u}(|\hat{y}_{s}|^2+|\hat{z}_{s}|^2)e^{\beta s}ds\right].
\end{eqnarray*}
Thus, if there exists a constant $\beta\geq2$ such that $\frac{2 K_1}{\beta}e^{\beta (T+l)}<1$, then there exists a constant $0<\gamma<1,$ such that
$$E\left[\int_{-l}^{T+u}(|\hat{Y}_s|^2+|\hat{Z}_s|^2)e^{\beta s} ds\right]\leq \gamma E\left[\int_{-l}^{T+u}(|\hat{y}_s|^2
+|\hat{z}_s|^2)e^{\beta s}ds\right].$$ Thus, $\|(\hat{Y}_s,\hat{Z}_s)\|_\beta<\sqrt{\gamma}\|(\hat{y}_s,\hat{z}_s)\|_\beta.$ Then by contraction mapping principle,  we can obtain (ii). Moreover, by Lemma 2.1, we have $(Y_t)_{t\in[0,T]}\in{\mathcal{S}}^2_{\cal{F}}(0,T)$. The proof is complete. \ $\Box$\\\\
\textbf{Remark 2} (i) From Theorem 2.2, it follows that if $g$ satisfies (A1) and (A3), BSDE (1) has a unique solution for a sufficiently small $l$ or a sufficiently small $K$, and if $g$ satisfies (A2) and (A3), BSDE (1) has a unique solution for a sufficiently small $T+l$ or a sufficiently small $K_1$,  Note that the existence and uniqueness of solution may be not true for arbitrary $K$ and $l$ (see Delong and Imkeller [6, Example 3.1]).

(ii) For simplicity, Theorem 2.2 is given for one-dimensional BSDE (1). In fact, since the proof of Theorem 2.2 is based on the estimate (Lemma 2.1) and fixed point theorem which both hold true in multidimensional case, we can know Theorem 2.2 also holds true for multidimensional BSDE (1). Theorem 2.2 generalizes the corresponding results in Delong and Imkeller [6], Peng and Yang [9] and Yang and Elliott [10] in some way. Cheridito and Nam [2] introduced a general BSDE whose coefficient depends on the solution on interval $[0,T]$. Our result in Theorem 2.2(ii) is for a BSDE whose coefficient depends on the solution on interval $[-l,T+u]$.

(iii) Since we use an estimate different from [6], the existence and uniqueness solution of BSDE (1) we obtain in Theorem 2.2(i) is independent of time horizon $T$, which is different from Delong and Imkeller [6, Theorem 2.1]. Such phenomena had also been discovered in Delong and Imkeller [6, Theorem 2.2] in a special case, and Cordoni et al. [3] for forward BSDEs with time-delayed coefficients.\\

By (i) in Remark 2, we can know if $g$ satisfies (A2) and (A3), BSDE (1) has a unique solution for a sufficiently small $T+l$. But the following example shows that if $g$ does not satisfy (A2), BSDE (1) may have no solution for any $T+l$. \\\\
\textbf{Example 2.3} Let $(\xi_t,\eta_t)$ satisfy terminal condition, $\delta\in[0,T].$ We consider
$$\left\{
    \begin{array}{ll}
      Y_t=\xi_T+\int_t^T \frac{1}{T-\delta}E[Y_{\delta}|{\cal{F}}_s])ds-\int_t^TZ_sdB_s,\ \ \ t\in[0,T);\\
      (Y_t, Z_t)=(\xi_t, \eta_t),\ t\in[T,T+u]\ \ \textmd{and} \ \
      (Y_t, Z_t)=(Y_0, 0), \ t\in[-l,0).
    \end{array}
  \right.$$
We can check that the coefficient of above BSDE does not satisfy (A2). If the above BSDE has solution, then we have $E[Y_\delta]=E[\xi_T]+E[Y_\delta]$. Thus $E[\xi_T]=0.$ Clearly, if $E[\xi_T]\neq0$, the BSDE will have no solution.\\

We give a simple example of BSDEs with time-advanced and -delayed coefficients.\\\\
\textbf{Example 2.4 } Let the time-delayed parameter $l=T$ and Lipschitz constant $K<\sqrt{\frac{1}{2e^{2T}}}.$  We consider
$$\left\{
    \begin{array}{ll}
      Y_t=B^2_T-T+1 +\int_t^T K(Y_{s-T}+E[Z_{s+u}|{\cal{F}}_s]-Z_s)ds-\int_t^TZ_sdB_s,\ \ \ t\in[0,T);\\
      (Y_t, Z_t)=(B^2_t-t+1, 2B_t),\ t\in[T,T+u]\ \ \textmd{and} \ \
      (Y_t, Z_t)=(Y_0, 0), \ t\in[-T,0).
    \end{array}
  \right.\eqno(7)$$
We can check that the coefficient of above BSDE satisfies (A1) and (A3). By Theorem 2.2(i), BSDE (7) has a unique solution. Now, we firstly consider the following BSDE (see [6])
$$\left\{
    \begin{array}{ll}
      Y_t^1=1+\int_t^TKY_{s-T}^1ds-\int_t^TZ_s^1dB_s,\ \ t\in[0,T);\\
      (Y_t, Z_t)=(1, 0),\ t\in[T,T+u]\ \ \textmd{and} \ \
      (Y_t, Z_t)=(Y_0, 0), \ t\in[-T,0).
    \end{array}
  \right.\eqno(8)$$
Since $Y^1_{s-T}=Y^1_0,\ s\in[0,T],$ we can check that
 $$(Y_t^1,Z_s^1)=\left\{\begin{array}{ll}(\frac{1-tK}{1-TK},0),\ \ t\in[0,T),\\
 (Y_t, Z_t)=(1, 0),\ t\in[T,T+u]\ \ \textmd{and} \ \
      (Y_t, Z_t)=(Y_0, 0), \ t\in[-T,0),
 \end{array}
 \right.$$ is a solution of BSDE (8). Now, we consider the following BSDE (see [9])
$$\left\{
    \begin{array}{ll}
      Y_t^2=B^2_T-T+\int_t^TK(E[Z_{s+u}^2|{\cal{F}}_s]-Z_s^2)ds-\int_t^TZ_s^2dB_s,\ \ t\in[0,T);\\
      (Y_t, Z_t)=(B^2_t-t, 2B_t),\ t\in[T,T+u]\ \ \textmd{and} \ \
      (Y_t, Z_t)=(Y_0, 0), \ t\in[-T,0).
    \end{array}
  \right.\eqno(9)$$
We can check
  $$(Y_t^1,Z_s^1)=\left\{\begin{array}{ll}(B^2_t-t, 2B_t),\ \ t\in[0,T),\\
 (Y_t, Z_t)=(B^2_t-t, 2B_t),\ t\in[T,T+u]\ \ \textmd{and} \ \
      (Y_t, Z_t)=(Y_0, 0), \ t\in[-T,0),
 \end{array}
 \right.$$
is a solution of BSDE (9). Set $(Y_t,Z_t):=(Y_t^1,Z_t^1)+(Y_t^2,Z_t^2), t\in[-T,T+u],$ then we can check that $(Y_t,Z_t)$ is the unique solution of BSDE (7).\\

The following is a continuous dependence
property of BSDEs (1).\\\\
\textbf{Proposition 2.5}\textit{\ Let $(\xi_t,\eta_t), (\xi'_t,\eta'_t)$ satisfy terminal condition, $\varphi_t, \varphi'_t\in L^2_{\cal{F}}(0,T;{\mathbf{R}})$ and $g$ satisfy (A1) and (A3). Let $f(t,\cdot,\cdot):=g(t,\cdot,\cdot)+\varphi_t$ and $f'(t,\cdot,\cdot):=g(t,\cdot,\cdot)+\varphi'_t.$ Suppose that there exists a constant $\beta\geq2$ such that $\frac{8K^2e^{\beta l}}{\beta}\leq 1,$ where $K$ and $l$ are the constants in (A1). Let $(\textbf{Y},\textbf{Z})$ and $(\textbf{Y}',\textbf{Z}')$ be the solutions of BSDE (1) with coefficients $f$ and $f',$ respectively, then for each $t\in[0,T],$ we have}
\begin{eqnarray*}
&&|Y_t-Y'_t|^2+E\left[\int_{t}^T(|Y_s-Y'_s|^2+|Z_s-Z'_s|^2)ds|{\cal{F}}_t\right]\\
&\leq&CE\left[|\xi_T-\xi'_T|^2+\int_{t}^T|\varphi_s-\varphi'_s|^2ds+
\int_{T}^{T+u}(|\xi_s-\xi'_s|^2+|\eta_s-\eta'_s)ds|{\cal{F}}_t\right]\\&&+C\int_{t-l}^{t}(|Y_s-Y'_s|^2+|Z_s-Z'_s|^2)ds,
\end{eqnarray*}
\textit{where $C>0$ is a constant depending only on $K$ and $l$.}\\\\
\textit{Proof.} Set $\hat{\textbf{Y}}=\textbf{Y}-\textbf{Y}',\ \hat{\textbf{Z}}=\textbf{Z}-\textbf{Z}',\ \hat{\xi}=\xi-\xi'$, $\hat{\eta}=\eta-\eta'$ and $\hat{g}=g(s,{\textbf{Y}},{\textbf{Z}})-g(s,{\textbf{Y}}',{\textbf{Z}}')$. By (4) and a similar argument as (6), we can deduce
\begin{eqnarray*}
&&|\hat{Y}_t|^2+E\left[\int_{t}^T(\frac{\beta}{2}|\hat{Y}_s|^2+|\hat{Z}_s|^2)e^{\beta (s-t)} ds|{\cal{F}}_t\right]\\
&\leq&\left[|\xi_T-\xi'_T|^2e^{\beta T}|{\cal{F}}_t\right]+\frac{4}{\beta}E\left[\int_{t}^T|\varphi_s-\varphi'_s|^2e^{\beta (s-t)}ds|{\cal{F}}_t\right]+\frac{4}{\beta}E\left[\int_{t}^T|\hat{g}|^2e^{\beta (s-t)}ds|{\cal{F}}_t\right]\\
&\leq&E\left[|\xi_T-\xi'_T|^2e^{\beta T}|{\cal{F}}_t\right]+\frac{4}{\beta}E\left[\int_{t}^T|\varphi_s-\varphi'_s|^2e^{\beta (s-t)}ds|{\cal{F}}_t\right]\\
&&+\frac{8K^2e^{\beta l}}{\beta}E\left[\int_{t-l}^{T+u}(|\hat{Y}_{s}|^2+|\hat{Z}_{s}|^2)e^{\beta (s-t)}ds|{\cal{F}}_t\right].
\end{eqnarray*}
Then if there exists a constant $\beta\geq2$ such that $\frac{8 K^2e^{\beta l}}{\beta}\leq1,$ we can finish this proof from the above inequality. \ $\Box$\\

In general, the comparison theorem of BSDE (1) may not true (see Delong and Imkeller [6] or Peng and Yang [9]). But under some restrict conditions, we have the following Proposition 2.5. \\\\
\textbf{Proposition 2.6}\textit{\ Let $(\xi_t,\eta_t), (\xi'_t,\eta'_t)$ satisfy terminal condition, and for each $t\in[0,T]$, $({\textbf{y}},{\textbf{z}})\in \overline{L}^2_{\cal{F}}(-l,T+u)\times\overline{L}^{2}_{\cal{F}}(-l,T+u;{\mathbf{R}}^d)$, $g(t,{\textbf{y}},{\textbf{z}}), g'(t,{\textbf{y}},{\textbf{z}})$ are both dependent only on $(t,\textbf{y},z_t)$ (we denote them by $g(t,{\textbf{y}},z_t)$ and $g'(t,{\textbf{y}},z_t)$, respectively) and both satisfy (A1)' and (A3). Suppose that}

\textit{(i) $\xi_t\geq \xi'_t$ for each $t\in[T,T+u];$}

\textit{(ii)  for each $t\in[0,T]$, $({\textbf{y}},{\textbf{z}})\in \overline{L}^2_{\cal{F}}(-l,T+u)\times\overline{L}^{2}_{\cal{F}}(-l,T+u;{\mathbf{R}}^d)$, we have $g(t,{\textbf{y}},z_t)\geq g'(t,{\textbf{y}},z_t)$, and for each $t\in[0,T]$, ${\textbf{y}},{\textbf{y}}'\in \overline{L}^2_{\cal{F}}(-l,T+u)$, we have $g'(t,{\textbf{y}},z_t)\geq g'(t,{\textbf{y}}',z_t),$ if $y_r\geq y'_r$ for each $r\in[0,T]$.}\\
\textit{Let $(\textbf{Y},\textbf{Z})$ and $(\textbf{Y}',\textbf{Z}')$ be the solutions of BSDEs (1) with parameters $(g,\xi,\eta)$ and $(g',\xi',\eta'),$ respectively. Then for a sufficiently small $l$ or a sufficiently small $K'$ given in (A1)', we have, for each $t\in[-l,T+u],$ $Y_t\geq Y'_t.$}\\\\
\textit{Proof.} The proof uses an iteration method from Peng and Yang [9]. Clearly, $g'(s,\textbf{Y},\cdot)$ satisfies (A1)'' (standard Lipschitz condition), then by Peng [8], the following BSDE
$$\left\{
    \begin{array}{ll}
      Y_t^1=\xi_T +\int_t^Tg'(s,\textbf{Y},Z^1_t)
             ds-\int_t^TZ_s^1dB_s,\ \ \ t\in[0,T);\\
      (Y_t^1, Z_t^1)=(\xi'_t, \eta'_t),\ t\in[T,T+u]\ \ \textmd{and} \ \
      (Y_t^1, Z_t^1)=(Y_0, 0), \ t\in[-l,0).
    \end{array}
  \right.$$
has a unique solution $(\textbf{Y}^1,\textbf{Z}^1),$ and by (ii) and comparison theorem (see Peng [8]), we further have $Y_t\geq Y^1_t.$ By Peng [8] again, the following BSDE
$$\left\{
    \begin{array}{ll}
      Y_t^2=\xi_T +\int_t^Tg'(s,\textbf{Y}^1,Z^2_t)
             ds-\int_t^TZ_s^2dB_s,\ \ \ t\in[0,T);\\
      (Y_t^2, Z_t^2)=(\xi'_t, \eta'_t),\ t\in[T,T+u]\ \ \textmd{and} \ \
      (Y_t^2, Z_t^2)=(Y_0, 0), \ t\in[-l,0).
    \end{array}
  \right.$$
has a unique solution $(\textbf{Y}^2,\textbf{Z}^2).$ Since $Y_t\geq Y^1_t,$ by (ii) and comparison theorem (see Peng [8]) again, we have $Y^1_t\geq Y^2_t.$ Similarly, for $n>2,$ the BSDE (1) with parameter $(g'(s,{\textbf{Y}}^{n-1},\cdot),\xi',\eta',l,u)$ has a unique solution $(\textbf{Y}^n,\textbf{Z}^n)$ and $Y^{n-1}_t\geq Y^n_t.$ Set $\hat{\textbf{Y}}^n:=\textbf{Y}^n-\textbf{Y}^{n-1},\ \hat{\textbf{Z}}^n:=\textbf{Z}^n-\textbf{Z}^{n-1}, n\geq1,$ and $\textbf{Y}^0:=\textbf{Y}$, $\textbf{Z}^0:=\textbf{Z}$. By the similar treatment as (6), we have
\begin{eqnarray*}
&&E\left[\int_{-l}^{T+u}(\frac{\beta}{2}|\hat{Y}_s^n|^2+|\hat{Z}_s^n|^2)e^{\beta s} ds\right]\\&\leq& \frac{2}{\beta}E\left[\int_0^T|g'(s,{\textbf{Y}}^{n-1},Z_t^{n})-g'(s,{\textbf{Y}}^{n-2},Z_t^{n-1})|^2e^{\beta s}ds\right]\\&\leq& \frac{4}{\beta}E\left[K'^2\int_0^T\left(\int_{-l}^{u}|\hat{Y}^{n-1}_{s+r}|^2\lambda'(dr)+|\hat{Z}_s^n|^2\right) e^{\beta s}ds\right]
\\&=&\frac{4}{\beta}E\left[K'^2\left(\int_{-l}^{u}\int_0^T|\hat{Y}^{n-1}_{s+r}|^2e^{\beta s}ds\lambda'(dr)+\int_0^T|\hat{Z}_s^n|^2 e^{\beta s}ds\right)\right]
\\&=&\frac{4}{\beta}E\left[K'^2\left(\int_{-l}^{u}e^{-r}\int_r^{T+r}|\hat{Y}^{n-1}_{v}|^2e^{\beta v}dv\lambda'(dr)+\int_0^T|\hat{Z}_s^n|^2 e^{\beta s}ds\right)\right]
\\&\leq& \frac{4K'^2}{\beta}e^{\beta l}E\left[\int_{-l}^{T+u}(|\hat{Y}^{n-1}_s|^2+|\hat{Z}_s^n|^2)e^{\beta s}ds\right], \ n\geq1.
\end{eqnarray*}
If $K'$ or $l$ is small enough such that there exists a constant $\beta\geq2$ such that $\frac{4K'^2e^{\beta l}}{\beta}<\frac{1}{3},$ then by the above inequality, there exists a constant $0<\gamma<1,$ such that
\begin{eqnarray*}
E\left[\int_{-l}^{T+u}(|\hat{Y}_s^n|^2+|\hat{Z}_s^n|^2)e^{\beta s} ds\right]&\leq&\gamma E\left[\int_{-l}^{T+u}|\hat{Y}^{n-1}_s|^2e^{\beta s} ds\right]\\
&\leq&\gamma^2 E\left[\int_{-l}^{T+u}|\hat{Y}_s^{n-2}|^2e^{\beta s} ds\right]\\&\cdot\cdot\cdot&\\&\leq&\gamma^{n-1}E\left[\int_{-l}^{T+u}|\hat{Y}_s^1|^2e^{\beta s} ds\right].
\end{eqnarray*}
By this inequality, we can get that $(\textbf{Y}^n, \textbf{Z}^n)_{n\geq1}$ is a Cauchy sequence in $\overline{L}^2_{\cal{F}}(-l,T+u)\times \overline{L}^{2}_{\cal{F}}(-l,T+u;{\mathbf{R}}^d).$ We denote the limit of $(\textbf{Y}^n, \textbf{Z}^n)_{n\geq1}$ in $\overline{L}^2_{\cal{F}}(-l,T+u)\times \overline{L}^{2}_{\cal{F}}(-l,T+u;{\mathbf{R}}^d)$ by $(\textbf{Y}'', \textbf{Z}'')$. Then by (A1)' and the similar treatment as (6), for each $t\in[0,T],$ we have
\begin{eqnarray*}
&&E\left[\int_t^Tg'(s,{\textbf{Y}}^{n-1},{{Z}}^{n}_s)ds-\int_t^Tg'(s,{\textbf{Y}}'',{{Z}}_s'')ds\right]^2\\&\leq& TE\left[\int_t^T|g'(s,{\textbf{Y}}^{n-1},{{Z}}^{n}_s)-g'(s,{\textbf{Y}}'',{{Z}}_s'')|^2ds\right]
\\&\leq& \frac{2 T K'^2}{\beta}e^{\beta l}E\left[\int_{t-l}^{T+u}(|{{Y}}^{n-1}_v-{{Y}}''_v|^2+|{{Z}}^{n}_v-{{Z}}''_v|^2)dv\right]\\
&\rightarrow&0,
\end{eqnarray*}
as $n\rightarrow\infty.$ We also have
\begin{eqnarray*}
E\left[\int_t^T{{Z}}^n_sdB_s-\int_t^T{{Z}}''_sdB_s\right]^2&=& E\left[\int_t^T|{{Z}}^n_s-{{Z}}''_s|^2ds\right]\\
&\rightarrow&0,
\end{eqnarray*}
as $n\rightarrow\infty.$ Then from the above two equations, we can deduce that for almost every $t\in[0,T]$, we have ${Y}''_t={Y}'_t.$ In view of $\textbf{Y}^n\searrow \textbf{Y}''$ and $Y_t\geq Y^1_t,$ we have for almost every $t\in[0,T]$, we have ${Y}_t\geq {Y}'_t.$ By the continuity of $(Y_t)_{\in[0,T]}$ and $(Y'_t)_{\in[0,T]},$ we have for each $t\in[0,T]$, $Y_t\geq Y'_t.$  The proof is complete. $\Box$
\section{SDEs with time-advanced and -delayed coefficients}
In this section, let $t_0\in[0,T]$ be a given number. We consider the functions $b$ and $\sigma:$
$${b}\left( \omega,t,x\right) : \Omega \times [0,T]\times \overline{L}^2_{\cal{F}}(t_0-l,T+u)\longmapsto {\mathbf{R}},\
{\textmd{and}} \ {\sigma}\left( \omega ,t,x\right) : \Omega \times [0,T]\times \overline{L}^2_{\cal{F}}(t_0-l,T+u)\longmapsto {\mathbf{R}}^d,$$ such that for each $\textbf{x}\in \overline{L}^2_{\cal{F}}(t_0-l,T+u),$
$({b}(t,\textbf{x}))_{t\in [0,T]}$ and $({\sigma}(t,\textbf{x}))_{t\in [0,T]}$ are both progressively measurable processes. We make the following assumptions:
\begin{itemize}
  \item (B1). There exist a constant $K_2>0$ and a probability measure $\lambda_1$ defined on $[-l,u]$, such that for each $t\in[0,T]$ and each $\textbf{x},\tilde{\textbf{x}}\in \overline{L}^2_{\cal{F}}(t_0-l,T+u),$
  \begin{eqnarray*}
  |{b}(t,\textbf{x})-{b}(t,\tilde{\textbf{x}}) |+|{\sigma}(t,\textbf{x})-{\sigma}(t,\tilde{\textbf{x}}) |
   \leq K_2E\left[\int_{-l}^{u}|x_{t+r}
    -\tilde{x}_{t+r}|\lambda_1(dr)|{\cal{F}}_t\right];
\end{eqnarray*}
  \item (B2). There exists a constant $K_3>0$, such that for each $\textbf{x},\tilde{\textbf{x}}\in \overline{L}^2_{\cal{F}}(t_0-l,T+u),$
  \begin{eqnarray*}
&&E\left[\int_{t_0-l}^{T+u}|{b}(t,\textbf{x})-{b}(t,\tilde{\textbf{x}}) |^2+|{\sigma}(t,\textbf{x})-{\sigma}(t,\tilde{\textbf{x}}) |^2dt\right]
   \leq K_3E\left[\int_{t_0-l}^{T+u}|x_{t}
    -\tilde{x}_{t}|^2dt\right];
\end{eqnarray*}
  \item (B3). For each $\textbf{x}\in \overline{L}^2_{\cal{F}}(t_0-l,T+u),$ $b(t,\textbf{x})\in L^2_{\cal{F}}(0,T)$ and $\sigma(t,\textbf{x})\in L^2_{\cal{F}}(0,T;{\mathbf{R}}^d).$
\end{itemize}
\textbf{Remark 3}\ \ Similar to (ii) in Remark 1, we can show (B1) implies (B2).\\

For convenience, if $\theta_t\in \overline{L}^2_{\cal{F}}(t_0-l,t_0)$ with $\theta_{t_0}\in{L}^2({\cal{F}}_{t_0}),$ we call $\theta_t$ satisfy initial condition.  Let $\theta_t$ satisfy initial condition. We consider the following SDE:
$$\left\{
    \begin{array}{ll}
      X_t=X_{t_0} +\int_{t_0}^tb(s,\textbf{X})
             ds+\int_{t_0}^t\sigma(s,\textbf{X}) dB_s,\ \ \ t\in(t_0,T];\\
      X_t=\theta_t,\ t\in[t_0-l,t_0]\ \ \textmd{and} \ \ X_t=X_T,\ t\in(T,T+u].
    \end{array}
  \right.\eqno(10)$$
where we denote $\{X_t\}_{t_0-l\leq t\leq T+u}$ by $\textbf{X}.$ The solution of SDE (10) is $\textbf{X}$ satisfying (10) and $(X_t)_{t\in[t_0,T]}\in {\cal{S}}^2_{\cal{F}}(t_0,T)$, which depends only on parameter $(b,\sigma, T, t_0, \theta, l,u).$ Clearly, the coefficients not only depend on the value of its solutions of the present but also the past and the future.\\

Motivated by Peng [8, Lemma 3.1], we have the following Lemma 3.1.\\\\
\textbf{Lemma 3.1}\textit{\  Let $\theta_t$ satisfy initial condition. Suppose that $b_0(s), \sigma_0(s)\in {L}^{2}_{\cal{F}}(0,T).$  Then the SDE with coefficients $b_0(s)$ and $\sigma_0(s)$ has a unique solution $X,$ and the following estimate}
\begin{eqnarray*}
\ \ \ \ \ \ \ \ \ \ \ \ \ \ \ \ &&E[|X_T|^2e^{-\beta (T+u)}|{\cal{F}}_{t_0}]+\frac{\beta}{2}E\left[\int_{t_0}^{T+u}|X_s|^2e^{-\beta s} ds|{\cal{F}}_{t_0}\right]\\&\leq& |X_{t_0}|^2e^{-\beta t_0}+E\left[\int_{t_0}^{T}(\frac{2}{\beta}|b_0(s)|^2+|\sigma_0(s)|^2)e^{-\beta s}ds|{\cal{F}}_{t_0}\right],\ \ \ \ \ \ \ \ \ \ \ \  \ \ \ \ \ \ \ \ \ \ \ (11)
\end{eqnarray*}
\emph{holds true for arbitrary constant $\beta>0$. We also have}
$$E\left[\sup_{{t_0}\leq t \leq T}|X_t|^2\right]\leq CE\left[|X_{t_0}|^2+\int_{t_0}^T(|b_0(s)|^2+|\sigma_0(s)|^2)ds\right],\eqno(12)$$
\emph{where $C>0$ is a constant depending only on $T$.}\\\\
\emph{Proof.} By the classic SDEs theory, this SDE has a unique solution. For an arbitrary constant $\beta>0,$ applying It\^{o}'s formula for $|X_s|^2e^{-\beta s}$ for $s\in[t_0,t],$ we can deduce for $t\in[t_0,T],$
\begin{eqnarray*}
&&|X_{t_0}|^2e^{-\beta t_0}+E\left[\int_{t_0}^{t}(-\beta|X_s|^2+|\sigma_0(s)|^2)e^{-\beta s}ds|{\cal{F}}_{t_0}\right]\\&=& E[|X_t|^2e^{-\beta t}|{\cal{F}}_{t_0}]-2E\left[\int_{t_0}^{t}X_sb_0(s)e^{-\beta s} ds|{\cal{F}}_{t_0}\right]\\
&\geq&E\left[|X_t|^2e^{-\beta t}|{\cal{F}}_{t_0}\right]-\frac{\beta}{2}E\left[\int_{t_0}^{t}|X_s|^2e^{-\beta s} ds|{\cal{F}}_{t_0}\right]-\frac{2}{\beta}E\left[\int_{t_0}^{t}|b_0(s)|^2e^{-\beta s} ds|{\cal{F}}_{t_0}\right].
\end{eqnarray*}
Then we have for $t\in[t_0,T],$
$$E[|X_t|^2e^{-\beta t}|{\cal{F}}_{t_0}]+\frac{\beta}{2}E\left[\int_{t_0}^{t}|X_s|^2e^{-\beta s} ds|{\cal{F}}_{t_0}\right]\leq |X_{t_0}|^2e^{-\beta t_0}+E\left[\int_{t_0}^{t}(\frac{2}{\beta}|b_0(s)|^2+|\sigma_0(s)|^2)e^{-\beta s}ds|{\cal{F}}_{t_0}\right],\eqno(13)$$
Since $X_t=X_T, t\in[T,T+u],$ we have
$$|X_T|^2e^{-\beta (T+u)}=|X_T|^2e^{-\beta T}-\beta\int_T^{T+u}|X_s|^2e^{-\beta s} ds.$$
From this and (13), we can get (11). By BDG inequality and (11), we can get (12). $\Box$\\

The following is the main result of this section. It gives two sufficient conditions, under which SDE (10) has a unique solution.\\\\
\textbf{Theorem 3.2}\textit{\  Let $\theta_t$ satisfy initial condition. Then we have}

\textit{(i) If $b$ and $\sigma$ both satisfy (B1) and (B3), and there exists a constant $\beta>0$ such that $\frac{4 K_2^2 e^{\beta u}}{\beta}(1+ \frac{2}{\beta})<1,$  then SDE (10) has a unique solution, where $K_2$ and $u$ are constants in (B1).}

\textit{(ii) If $b$ and $\sigma$ both satisfy (B2) and (B3), and there exists a constant $\beta>0$ such that $\frac{2 K_3 e^{\beta(T+u-t_0)}}{\beta}\left(1+ \frac{2}{\beta}\right) <1$,  then SDE (10) has a unique solution, where $K_3$ and $l$ are constants in (B2).}\\\\
\emph{Proof.} We still can prove this theorem using a contractive method like Theorem 2.2. Clearly, $\overline{{L}}^2_{\cal{F}}(t_0-l,T+u)$ is a Banach space with the norm:
$$\|\cdot\|_{-\beta}=E\left[\int_{t_0-l}^{T+u}|\cdot|^2e^{-\beta s} ds\right]^{\frac{1}{2}},$$
where $\beta>0$ is a constant. By (B3) and Lemma 3.1, we can define a mapping $\phi$ from
$\overline{{L}}^2_{\cal{F}}(t_0-l,T+u)$
into itself by setting $\textbf{X}:=\phi(\textbf{x}),$
where $\textbf{X}$ is the solution of the SDE
$$\left\{
    \begin{array}{ll}
      X_t=X_{t_0} +\int_{t_0}^tb(s,\textbf{x})
             ds+\int_{t_0}^t\sigma(s,\textbf{x}) dB_s,\ \ \ t\in(t_0,T];\\
      X_t=\theta_t,\ t\in[t_0-l,t_0]\ \ \textmd{and} \ \ X_t=X_T,\ t\in(T,T+u].
    \end{array}
  \right.$$
For $\textbf{x}^1, \textbf{x}^2\in\overline{{L}}^2_{\cal{F}}(t_0-l,T+u),$ let $\textbf{X}^i:=\phi(\textbf{x}^i),\ \  i=1,2.$ We set
$\hat{X}_t:={X}^1_t-{X}^2_t$ and $\hat{x}_t:={x}^1_t-{x}^2_t.$

\textbf{Proof of (i):}

By Lemma 3.1, (B1) and Fubini's theorem, we have
\begin{eqnarray*}
&&E\left[\int_{t_0}^{T+u}\frac{\beta}{2}|\hat{X}_s|^2e^{-\beta s} ds\right]\\ &\leq& E\left[\int_{t_0}^T (\frac{2}{\beta}|b(s,\textbf{x}^1)-b(s,\textbf{x}^2)|^2+|\sigma(s,\textbf{x}^1)-\sigma(s,\textbf{x}^2)|^2)e^{-\beta s}ds\right]\\&\leq&2K_2^2\left(1+ \frac{2}{\beta}\right)E\left[\int_{t_0}^T\int_{-l}^u|\hat{x}_{s+r}|^2\lambda_1(dr)e^{-\beta s}ds\right]
\\&\leq&2K_2^2\left(1+ \frac{2}{\beta}\right)E\left[\int_{-l}^ue^{\beta r}\int_{t_0}^T|\hat{x}_{s+r}|^2 e^{-\beta (s+r)}ds\lambda_1(dr)\right]
\\&\leq&2K_2^2\left(1+ \frac{2}{\beta}\right)e^{\beta u}E\left[\int_{-l}^u\int_{t_0+r}^{T+r}|\hat{x}_{v}|^2e^{-\beta v}dv\lambda_1(dr)\right]\\&\leq&2K_2^2 e^{\beta u}\left(1+ \frac{2}{\beta}\right)E\left[ \int_{t_0-l}^{T+u}|\hat{x}_{v}|^2e^{-\beta v}dv\right].
\ \ \ \ \ \  \ \ \ \ \ \ \ \  \ \ \ \ \ \ \ \  \ \ \ \ \ \ \ \  \ \ \ \ \ \ \ \  \ \   \ \ \ \ \ \ \ \  \ \ \ \ \ \ \ \  \ \ (14)
\end{eqnarray*}
Thus we have
\begin{eqnarray*}
E\left[\int_{t_0-l}^{T+u}|\hat{X}_s|^2e^{-\beta s} ds\right]\leq\frac{4 K_2^2 e^{\beta u}}{\beta}\left(1+ \frac{2}{\beta}\right)E\left[\int_{t_0-l}^{T+u}|\hat{x}_{v}|^2e^{-\beta v}ds\right].
\end{eqnarray*}
Thus, if there exists a constant $\beta>0$ such that $\frac{4 K_2^2 e^{\beta u}}{\beta}\left(1+ \frac{2}{\beta}\right)<1$, then there exists a constant $0<\gamma<1,$ such that $\|\hat{X}_s\|_{-\beta}<\sqrt{\gamma}\|\hat{x}_s\|_{-\beta}.$ Then by contraction mapping principle, we can obtain (i). Moreover, by Lemma 3.1, we have $(X_t)_{t\in[t_0,T]}\in{\mathcal{S}}^2_{\cal{F}}(t_0,T)$.

\textbf{Proof of (ii):}

By Lemma 3.1 and (B2), we have
\begin{eqnarray*}
&&E\left[\int_{t_0}^{T+u}\frac{\beta}{2}|\hat{X}_s|^2e^{-\beta s} ds\right]\\ &\leq& E\left[\int_{t_0}^T (\frac{2}{\beta}|b(s,\textbf{x}^1)-b(s,\textbf{x}^2)|^2+|\sigma(s,\textbf{x}^1)-\sigma(s,\textbf{x}^2)|^2)e^{-\beta s}ds\right]\\&\leq&K_3\left(1+ \frac{2}{\beta}\right)e^{-\beta t_0}E\left[\int_{t_0-l}^{T+u}|\hat{x}_{s}|^2ds\right]\\&\leq&K_3\left(1+ \frac{2}{\beta}\right)e^{\beta(T+u- t_0)}E\left[\int_{t_0-l}^{T+u}|\hat{x}_{s}|^2e^{-\beta s}ds\right].
\end{eqnarray*}
Thus we have
\begin{eqnarray*}
E\left[\int_{t_0-l}^{T+u}|\hat{X}_s|^2e^{-\beta s} ds\right]\leq \frac{2K_3e^{\beta(T+u- t_0)}}{\beta}\left(1+ \frac{2}{\beta}\right)E\left[\int_{t_0-l}^{T+u}|\hat{x}_{v}|^2e^{-\beta v}ds\right].
\end{eqnarray*}
Thus, if there exists a constant $\beta>0$ such that $\frac{2K_3e^{\beta(T+u- t_0)}}{\beta}\left(1+ \frac{2}{\beta}\right)<1,$ then there exists a constant $0<\gamma<1,$ such that $\|\hat{X}_s\|_{-\beta}<\sqrt{\gamma}\|\hat{x}_s\|_{-\beta}.$ Then by contraction mapping principle, we can obtain (ii). Moreover, by Lemma 3.1, we have  $(X_t)_{t\in[t_0,T]}\in{\mathcal{S}}^2_{\cal{F}}(t_0,T)$. The proof is complete. \ $\Box$\\\\
\textbf{Remark 4} (i) From Theorem 3.2, it follows that if $b$ and $\sigma$ both satisfy (B1) and (B3), then SDE (10) has a unique solution for a sufficiently small $u$ or a sufficiently small $K_2,$ and if $b$ and $\sigma$ both satisfy (B2) and (B3), SDE (10) has a unique solution for a sufficiently small $T+u-t_0$ or a sufficiently small $K_3.$ But the following example shows SDE (10) may have no solution for some $K_2$ and $u.$ We consider the SDE $$\left\{
    \begin{array}{ll}
      X_t=X_{0} +\int_{0}^tK_2E[X_{s+u}|{\cal{F}}_s]ds+\int_{0}^tK_2\textbf{I}_ddB_s,\ \ \ t\in(0,T];\\
      X_0=a,\  \ \textmd{and} \ \ X_t=X_T,\ t\in(T,T+u].
    \end{array}
  \right.$$
  Let $u=T,$ we have $X_T=a+\int_{0}^TK_2E[X_T|{\cal{F}}_s]ds+\int_{0}^TK_2\textbf{I}_ddB_s,$ then we have $a=E[X_T](1-TK_2).$ If $TK_2=1$ and $a\neq0,$ the above SDE will have not solution.

(ii) For simplicity, Theorem 3.2 is given for one-dimensional SDE (10). In fact, since the proof of Theorem 3.2 is based on the estimate (Lemma 3.1) and fixed point theorem which both hold true in multidimensional case, we can know Theorem 3.2 also holds true for multidimensional SDE (10). Theorem 3.2 generalizes the classic result for time-delayed SDEs and the corresponding result in Chen and Huang [1].\\

By (i) in Remark 4, we can know if $b$ and $\sigma$ both satisfy (B2) and (B3), SDE (10) has a unique solution for a sufficiently small $T+u-t_0$. But the following example shows that if $b$ does not satisfy (B2), SDE (10) may have no solution for any $T+u-t_0$. \\\\
\textbf{Example 3.3} Given $\delta\in[0,T]$, we consider the SDE $$\left\{
    \begin{array}{ll}
      X_t=X_{t_0} +\int_{t_0}^t\frac{1}{\delta-t_0}E[X_{\delta}|{\cal{F}}_s]ds+\int_{t_0}^t\textbf{I}_ddB_s,\ \ \ t\in(t_0,T];\\
      X_t=a\ ,\ t\in[t_0-l,t_0]\ \ \textmd{and} \ \ X_t=X_T,\ t\in(T,T+u].
    \end{array}
  \right.$$
We can check the coefficient $b$ does not satisfy (B2). If the above SDE has solution, then we have $E[X_\delta]=a+E[X_\delta]$. Thus $a=0.$ Clearly, if $a\neq0$, the above SDE will have no solution.\\

By (13) and the similar arguments as (14), we can get the following continuous dependence
property of the SDEs (10). \\\\
\textbf{Proposition 3.4}\textit{\ Let $\theta_t, \theta'_t$ satisfy initial condition, $\varphi_t, \varphi'_t\in L^2_{\cal{F}}(0,T;{\mathbf{R}}),$ and $\phi_t, \phi'_t\in L^2_{\cal{F}}(0,T;{\mathbf{R}}^d).$ Set $\hat{b}(t,\cdot):=b(t,\cdot)+\varphi_t,\ \hat{b}'(t,\cdot):=b(t,\cdot)+\varphi'_t$ and $\hat{\sigma}(t,\cdot):=\sigma(t,\cdot)+\phi_t, \hat{\sigma}(t,\cdot):=\sigma(t,\cdot)+\phi'_t.$ Suppose that  $b$ and $\sigma$ both satisfy (B1) and (B3), and there exists a constant $\beta>0$ such that $\frac{8K_1^2 e^{\beta u}}{\beta}\left(1+ \frac{2}{\beta}\right)<1.$ If $\textbf{X}$ and $\textbf{X}'$ be the solutions of SDEs (11) with coefficients $(\hat{b},\hat{\sigma})$ and $(\hat{b}',\hat{\sigma}')$, respectively, then for each $t\in[0,T],$ we have}
\begin{eqnarray*}
 &&E[|X_t-X'_t|^2|{\cal{F}}_{t_0}]+E\left[\int_{t_0}^t|X_s-X'_s|^2ds|{\cal{F}}_{t_0}\right]\\
&\leq&CE\left[\int_{t_0}^t(|\varphi_s-\varphi'_s|^2+|\phi_s-\phi'_s|^2)ds+\int_{t}^{t+u}|X_s-X'_s|^2ds|{\cal{F}}_{t_0}\right]+
C\int_{t_0-l}^{t_0}|\theta_s-\theta'_s|^2ds
\\&&+C|\theta_{t_0}-\theta'_{t_0}|^2.
\end{eqnarray*}
\textit{where $C>0$ is a constant depending only on $K_2$ and $u$.}\\

\emph{Proof.} Set $\hat{X}_t=X_t-X'_t$, $\hat{\varphi}_t=\varphi_t-\varphi'_t$ and $\hat{\phi}_t=\phi_t-\phi'_t$. By (13) and the similar arguments as (14),  we can get
\begin{eqnarray*}
&&E[|\hat{X}_t|^2e^{-\beta t}|{\cal{F}}_{t_0}]+\frac{\beta}{2}E\left[\int_{t_0}^{t}|\hat{X}_s|^2e^{-\beta s} ds|{\cal{F}}_{t_0}\right]
\\ &\leq& |\hat{X}_{t_0}|^2e^{-\beta t_0}+E\left[\int_{t_0}^{t}(\frac{2}{\beta}|\hat{b}(s,\textbf{X})-\hat{b}'(s,\textbf{X})|^2+|\hat{\sigma}(s,\textbf{X})-\hat{\sigma}'(s,\textbf{X})|^2)e^{-\beta s}ds|{\cal{F}}_{t_0}\right]
\\ &\leq& |\hat{X}_{t_0}|^2e^{-\beta t_0}+2E\left[\int_{t_0}^{t}(\frac{2}{\beta}|\hat{\varphi}_s|^2+|\hat{\phi}_s|^2)e^{-\beta s}ds|{\cal{F}}_{t_0}\right]\\&&+2E\left[\int_{t_0}^{t}(\frac{2}{\beta}|{b}(s,\textbf{X})-{b}'(s,\textbf{X}')|^2+|{\sigma}(s,\textbf{X})-{\sigma}'(s,\textbf{X}')|^2)e^{-\beta s}ds|{\cal{F}}_{t_0}\right]
\\ &\leq& |\hat{X}_{t_0}|^2e^{-\beta t_0}+2E\left[\int_{t_0}^{t}(\frac{2}{\beta}|\hat{\varphi}_s|^2+|\hat{\phi}_s|^2)e^{-\beta s}ds|{\cal{F}}_{t_0}\right]\\&&+2\left(1+ \frac{2}{\beta}\right) E\left[\int_{t_0}^{t}(|{b}(s,\textbf{X})-{b}'(s,\textbf{X}')|^2+|{\sigma}(s,\textbf{X})-{\sigma}'(s,\textbf{X}')|^2)e^{-\beta s}ds|{\cal{F}}_{t_0}\right]
\\ &\leq& |\hat{X}_{t_0}|^2e^{-\beta t_0}+2E\left[\int_{t_0}^{t}(\frac{2}{\beta}|\hat{\varphi}_s|^2+|\hat{\phi}_s|^2)e^{-\beta s}ds|{\cal{F}}_{t_0}\right]\\&&+4K_2^2\left(1+ \frac{2}{\beta}\right)E\left[\int_{t_0}^{t}\int_{-l}^u(|\hat{X}_{s+r}|^2)\lambda_1(dr)e^{-\beta s}ds|{\cal{F}}_{t_0}\right]
\\&\leq&|\hat{X}_{t_0}|^2e^{-\beta t_0}+2\left(1+ \frac{2}{\beta}\right)E\left[\int_{t_0}^{t}(|\hat{\varphi}_s|^2+|\hat{\phi}_s|^2)e^{-\beta s}ds|{\cal{F}}_{t_0}\right]\\&&+4K_2^2\left(1+ \frac{2}{\beta}\right)E\left[\int_{-l}^u e^{\beta r}\int_{t_0}^{t}(|\hat{X}_{s+r}|^2)e^{-\beta (s+r)}ds\lambda_1(dr)|{\cal{F}}_{t_0}\right]
\\&\leq&|\hat{X}_{t_0}|^2e^{-\beta t_0}+2\left(1+ \frac{2}{\beta}\right)E\left[\int_{t_0}^{t}(|\hat{\varphi}_s|^2+|\hat{\phi}_s|^2)e^{-\beta s}ds|{\cal{F}}_{t_0}\right]\\&&+4K_2^2\left(1+ \frac{2}{\beta}\right)E\left[\int_{-l}^u e^{\beta r}\int_{t_0+r}^{t+r}(|\hat{X}_{v}|^2)e^{-\beta v}dv\lambda_1(dr)|{\cal{F}}_{t_0}\right]
\\&\leq&|\hat{X}_{t_0}|^2e^{-\beta t_0}+2\left(1+ \frac{2}{\beta}\right)E\left[\int_{t_0}^{t}(|\hat{\varphi}_s|^2+|\hat{\phi}_s|^2)e^{-\beta s}ds|{\cal{F}}_{t_0}\right]\\&&+4K_2^2 e^{\beta u}\left(1+ \frac{2}{\beta}\right)E\left[\int_{t_0-l}^{t+u}|\hat{X}_{v}|^2e^{-\beta v}dv|{\cal{F}}_{t_0}\right].
\end{eqnarray*}
If there exists a constant $\beta>0$ such that $\frac{8K_2^2 e^{\beta u}}{\beta}\left(1+ \frac{2}{\beta}\right)<1,$ then by above inequality, we can complete this proof. \ \ $\Box$\\

In general, the comparison theorem of SDE (10) may not true (see Yang et al. [11, Example 3.2 and Example 3.3]). But we can get the following Proposition 3.5. \\\\
\textbf{Proposition 3.5}\textit{\ Let $\theta_t, \theta'_t$ satisfy initial condition, and $b(s,\cdot), b'(s,\cdot)$ satisfy (B1) and (B3), and for each $\textbf{x}\in \overline{L}^2_{\cal{F}}(-l,T+u),$ $b(t,\textbf{x})$ and $ b'(t,\textbf{x})$ are both continuous in $t$. Let $\sigma_t$ be a continuous process and $\sigma_t\in L^2_{\cal{F}}(0,T;{\mathbf{R}}^d).$ Suppose that}

\textit{(i) $\theta_t\geq \theta'_t$ for each $t\in[t_0-l,t_0];$}

\textit{(ii) for each $t\in[0,T]$, $\textbf{x}\in \overline{L}^2_{\cal{F}}(t_0-l,T+u),$ we have $b(t,\textbf{x})\geq b'(s,\textbf{x}),$ and for each $t\in[0,T]$, $\textbf{x},\textbf{x}'\in \overline{L}^2_{\cal{F}}(t_0-l,T+u),$ we have $b'(t,\textbf{x})\geq b'(t,\textbf{x}'),$ if $x_{r}\geq x'_{r}$ for each $r\in[t_0,T].$}\\
\textit{Let SDEs (10) with coefficients $(b,\sigma)$ and $(b',\sigma)$ have a unique solution $\textbf{X}$ and $\textbf{X}'$, respectively. Then for a sufficiently small $u$ or a sufficiently small $K_2$ given in (B1), we have, for each $t\in[t_0-l,T+u],$ $X_t\geq X'_t.$}\\\\
\textit{Proof.} Clearly, the following SDE
$$\left\{
    \begin{array}{ll}
      X_t^1=X_{t_0}^1 +\int_{t_0}^tb'(s,\textbf{X})
             ds+\int_{t_0}^t\sigma_sdB_s,\ \ \ t\in(t_0,T];\\
      X_t^1=\theta_t',\ t\in[t_0-l,t_0]\ \ \textmd{and} \ \ X_t^1=X_T^1,\ t\in(T,T+u].
    \end{array}
  \right.$$
has a unique solution $\textbf{X}^1.$ By (ii) and comparison theorem (see Yang et al. [11, Lemma 2.4]), we can get $X_t\geq X^1_t.$ Clearly, the following SDE
$$\left\{
    \begin{array}{ll}
      X_t^2=X_{t_0}^1 +\int_{t_0}^tb'(s,\textbf{X}^1)
             ds+\int_{t_0}^t\sigma_sdB_s,\ \ \ t\in(t_0,T];\\
      X_t^2=\theta_t',\ t\in[t_0-l,t_0]\ \ \textmd{and} \ \ X_t^2=X_T^2,\ t\in(T,T+u].
    \end{array}
  \right.$$
has a unique solution $\textbf{X}^2.$ Since $X_t\geq X^1_t,$ by (ii) and comparison theorem again, we have $X^1_t\geq X^2_t.$ Similarly, for $n>2,$ the SDE (10) with parameter  $(b'(s,\textbf{X}^{n-1}),\sigma_s, \theta',l,u)$ has a unique solution $\textbf{X}^n$ and $X^{n-1}_t\geq X^n_t.$ Set $X_t^0=X_t$ and $\hat{X}^n_t=X_t^n-X_t^{n-1},\ n\geq1$, and by (11) and the same treatment as (14), we can get
\begin{eqnarray*}
E\left[\int_{t_0-l}^{T+u}\frac{\beta}{2}|\hat{X}_s^n|^2e^{-\beta s} ds\right]&\leq& E\left[\int_{t_0}^T (\frac{2}{\beta}|b'(s,\textbf{X}^{n-1})-b'(s,\textbf{X}^{n-2})|^2)e^{-\beta s}ds\right]\\&\leq&K_2^2\frac{2}{\beta}E\left[\int_{t_0}^T\int_{-l}^u|\hat{X}_{s+r}^{n-1}|^2\lambda_1(dr)e^{-\beta s}ds\right]
\\&\leq&\frac{2 K_2^2 e^{\beta u}}{\beta}E\left[ \int_{t_0-l}^{T+u}|\hat{X}^{n-1}_{v}|^2e^{-\beta v}dv\right].
\end{eqnarray*}
If $K'$ or $l$ is small enough such that there exists a constant $\beta>0$ such that $\frac{4K_2^2 e^{\beta u}}{\beta^2}<1,$ and SDEs (10) with coefficients $(b,\sigma)$ and $(b',\sigma),$ both have a unique solution, then by the above inequality, there exists a constant $0<\gamma<1,$ such that
\begin{eqnarray*}
E\left[\int_{t_0-l}^{T+u}|\hat{X}_s^n|^2e^{\beta s} ds\right]&\leq&\gamma E\left[ \int_{t_0-l}^{T+u}|\hat{X}^{n-1}_{s}|^2e^{-\beta s}ds\right]
\\&\cdot\cdot\cdot&\\&\leq&\gamma^{n-1}E\left[\int_{-l}^{T+u}|\hat{X}_s^{1}|^2e^{\beta s} ds\right]
\end{eqnarray*}
By this inequality, we can get that $(\textbf{X}^n)_{n\geq1}$ is a Cauchy sequence in $\overline{L}^2_{\cal{F}}(t_0-l,T+u),$ and denote its limit in $\overline{L}^2_{\cal{F}}(t_0-l,T+u)$ by $\textbf{X}''$. Then for each $t\in[t_0,T+u],$ we have
\begin{eqnarray*}
E\left[\int_{t_0}^{t}b'(s,\textbf{X}^{n-1})ds-\int_{t_0}^{t}b'(s,\textbf{X}'')ds\right]^2&\leq&(t-t_0)E\left[\int_{t_0}^t (|b'(s,\textbf{X}^{n-1})-b'(s,\textbf{X}'')|^2)ds\right]\\&\leq&(t-t_0)K_2^2E\left[\int_{t_0}^t\int_{-l}^u|X_{s+r}^{n-1}-X_{s+r}''|^2\lambda_1(dr)ds\right]
\\&\leq&(t-t_0)K_2^2E\left[\int_{t_0-l}^{t+u}|X_{v}^{n-1}-X_{v}''|^2dv\right]
\\&\rightarrow&0
\end{eqnarray*}
as $n\rightarrow\infty$. From the above inequality, we can deduce that for almost every $t\in[0,T]$, $X_t''=X_t'$. In view of $\textbf{X}^n\searrow \textbf{X}''$ and $X_t\geq X^1_t,$ we have for almost every $t\in[0,T]$, $X_t\geq X'_t.$ By the continuity of $(X_t)_{t\in[t_0,T]}$ and $(X'_t)_{t\in[t_0,T]}$, we have for every $t\in[0,T]$, $X_t\geq X'_t.$ The proof is complete. $\Box$
\section{A duality between BSDEs and SDEs with time-advanced and -delayed coefficients}
In this section, we will give a duality between BSDEs and SDEs with time-advanced and -delayed coefficients. Let $(\xi_t,\eta_t)$ satisfy terminal condition, we consider linear BSDE
$$\left\{
    \begin{array}{ll}
      Y_t=\xi_T +\int_t^T(b_sY_s+\bar{b}_{s}Y_{s-l}+E[\underline{b}_{s+u}Y_{s+u}|{\cal{F}}_s]+\sigma_s Z_{s}+\bar{\sigma}_{s} Z_{s-l}+E[\underline{\sigma}_{s+u}Z_{s+u}|{\cal{F}}_s]+c_s)
             ds\\ \ \ \ \ \ \ -\int_t^TZ_s dB_s, \ t\in[0,T);\\
      (Y_t, Z_t)=(\xi_t, \eta_t),\ t\in[T,T+u]\ \ \textmd{and}\ \ (Y_t, Z_t)=(Y_0, 0), \ t\in[-l,0),
    \end{array}
  \right.\eqno(15)$$
and linear SDE
$$\left\{
    \begin{array}{ll}
      X_r=X_{t} +\int_{t}^r(b_sX_s+\underline{{b}}_sX_{s-u}+E[\bar{{b}}_{s+l}X_{s+l}|{\cal{F}}_s])
             ds\\
              \ \ \ \ \ \ \ +\int_{t}^r(\sigma_s X_s+\underline{{\sigma}}_sX_{s-u}+E[\bar{{\sigma}}_{s+l}X_{s+l}|{\cal{F}}_s])
             dB_s,\ \ \ r\in(t,T];\\
      X_{t}=1\ \ \textmd{and} \ \ X_r=0,\ r\in(T,T+l],\ \ \textmd{and}\ \ X_r=0,\ r\in[t-u,t),
    \end{array}
  \right.\eqno(16)$$
where $b_t,\bar{b}_t,\underline{b}_t, c_t\in{L}^2_{\cal{F}}(0,T+u+l)$ and $\sigma_t,\bar{\sigma}_t,\underline{\sigma}_t\in{L}^2_{\cal{F}}(0,T+u+l;{\textbf{R}}^d)$ are all uniformly bounded. Here, we let advanced time (rep. delayed time) in BSDE (15) equal to delayed time (rep. advanced time) in SDE (16), in order to obtain a duality between them.  Clearly, the coefficient of BSDE (15) satisfies (A1) and (A3), and the coefficients of SDE (16) satisfy (B1) and (B3). By Theorem 2.2 and Theorem 3.2, BSDE (15) and SDE (16) have a unique solution, respectively, for a sufficiently small $l.$\\\\
\textbf{Theorem 4.1}\textit{\ \ Suppose that BSDE (15) and SDE (16) both have a unique solution, respectively. Then we have }
$$Y_t=E\left[X_T\xi_T+\int_T^{T+u}X_{s-u}(\underline{b}_{s}\xi_s+\underline{\sigma}_{s}\eta_s)ds
+\int_t^{t+l}X_s(\bar{{b}}_{s}Y_{s-l}+\bar{{\sigma}}_{s}Z_{s-l})ds+\int_t^Tc_sX_sds|{\cal{F}}_{t}\right].
$$
\textit{In particular, if $E\int_0^{l}X_s\bar{{b}}_{s}ds\neq1,$ we have the following closed formula}
$$Y_0=\left(1-E\left[\int_0^{l}X_s\bar{{b}}_{s}ds\right]\right)^{-1}E\left[X_T\xi_T
+\int_T^{T+u}X_{s-u}(\underline{b}_{s}\xi_s+\underline{\sigma}_{s}\eta_s)ds+\int_0^Tc_sX_sds\right].\eqno(17)$$
\textit{Proof.} Applying It\^{o}'s formula to $X_sY_s$ for $s\in[t,T]$ and taking conditional expectation, we have
\begin{eqnarray*}
&&E[X_TY_T|{\cal{F}}_{t}]-X_tY_t
\\&=&E\left[\int_t^{T}(\underline{b}_{s}X_{s-u}Y_s-E\left[\underline{b}_{s+u}X_{s}Y_{s+u}|{\cal{F}}_{s}\right])ds
+\int_t^{T}(E\left[\bar{{b}}_{s+l}X_{s+l}Y_s|{\cal{F}}_{s}\right]-\bar{{b}}_{s}X_{s}Y_{s-l})ds|{\cal{F}}_{t}\right]\\
&&+E\left[\int_t^{T}(\underline{\sigma}_{s}X_{s-u}Z_s-E\left[\underline{\sigma}_{s+u}X_{s}Z_{s+u}|{\cal{F}}_{s}\right])ds
+\int_t^{T}(E\left[\bar{{\sigma}}_{s+l}X_{s+l}Z_s|{\cal{F}}_{s}\right]-\bar{{\sigma}}_{s}X_{s}Z_{s-l})ds|{\cal{F}}_{t}\right]\\
&&-E\left[\int_t^Tc_sX_sds|{\cal{F}}_{t}\right]
\\&=&E\left[\int_t^{T}(\underline{b}_{s}X_{s-u}Y_s-\underline{b}_{s+u}X_{s}Y_{s+u})ds
+\int_t^{T}(\bar{{b}}_{s+l}X_{s+l}Y_s-\bar{{b}}_{s}X_{s}Y_{s-l})ds|{\cal{F}}_{t}\right]\\
&&+E\left[\int_t^{T}(\underline{\sigma}_{s}X_{s-u}Z_s-\underline{\sigma}_{s+u}X_{s}Z_{s+u})ds
+\int_t^{T}(\bar{{\sigma}}_{s+l}X_{s+l}Z_s-\bar{{\sigma}}_{s}X_{s}Z_{s-l})ds|{\cal{F}}_{t}\right]\\
&&-E\left[\int_t^Tc_sX_sds|{\cal{F}}_{t}\right]
\\&=&E\left[\left(\int_t^{T}\underline{b}_{s}X_{s-u}Y_sds-\int_{t+u}^{T+u}\underline{b}_{s}X_{s-u}Y_{s}ds\right)
+\left(\int_{t+l}^{T+l}\bar{{b}}_{s}X_{s}Y_{s-l}ds-\int_t^{T}\bar{{b}}_{s}X_{s}Y_{s-l}ds\right)|{\cal{F}}_{t}\right]\\
&&+E\left[\left(\int_t^{T}\underline{\sigma}_{s}X_{s-u}Z_sds-\int_{t+u}^{T+u}\underline{\sigma}_{s}X_{s-u}Z_{s}ds\right)
+\left(\int_{t+l}^{T+l}\bar{{\sigma}}_{s}X_{s}Z_{s-l}-\int_t^{T}\bar{{\sigma}}_{s}X_{s}Z_{s-l}ds\right)|{\cal{F}}_{t}\right]\\
&&-E\left[\int_t^Tc_sX_sds|{\cal{F}}_{t}\right].
\end{eqnarray*}
Since $X_t=1, X_s=0$ for $s\in(T,T+l]$ and $X_s=0$ for $s\in[t-u,t),$ we can get
$$Y_t=E\left[X_T\xi_T+\int_T^{T+u}X_{s-u}(\underline{b}_{s}\xi_s+\underline{\sigma}_{s}\eta_s)ds
+\int_t^{t+l}X_s(\bar{{b}}_{s}Y_{s-l}+\bar{{\sigma}}_{s}Z_{s-l})ds+\int_t^Tc_sX_sds|{\cal{F}}_{t}\right].
$$
Since $Y_s=Y_0, Z_s=0$ for $s\in[-l,0),$ we have
$$Y_0=E\left[X_T\xi_T+\int_T^{T+u}X_{s-u}(\underline{b}_{s}\xi_s+\underline{\sigma}_{s}\eta_s)ds
+Y_0\int_0^{l}X_s\bar{{b}}_{s}ds+\int_0^Tc_sX_sds|{\cal{F}}_{t}\right].$$
From this, we can get (17). The proof is complete. \ $\Box$\\

\textbf{Acknowledgements.}\ \ The authors would like to thank the anonymous referee for valuable
comments and suggestions.

\end{document}